\documentclass[12pt]{article}
\usepackage{pstricks,pst-node,pst-text,pst-3d}

  \usepackage{amssymb}
  \usepackage{amsmath}
\usepackage{graphics}
\usepackage[dvips]{epsfig}
\usepackage{times}

\newtheorem{lemma}{Lemma}[section]

\newtheorem{theorem}{Theorem}[section]
\newtheorem{ex}{Example}[section]

\def\R{{\mathbb R}}

\def\Rem{\textsc{Remark}}
\def\Pf{{\it Proof.$\;$}}

\def\qed{\hfill $\diamond$}

\def\({\langle}
\def\){\rangle}
\def\mb{\boldsymbol}

\def\cC{{\mathcal C}}

\def\cF{{\mathcal F}}
\def\cG{{\mathcal G}}

\def\cN{{\mathcal N}}

\def\1{\mb1}

\def\v0{{\bf 0}}

\def\ov{\overline}

\begin{document}

\title{\bf\Large Least Square Approximations and\\ Linear Values of Cooperative Games}

\author{Ulrich FAIGLE\thanks{Mathematisches Institut, Universit\"at zu K\"oln,
    Weyertal 80, 50931 K\"oln, Germany. Email: \texttt{faigle@zpr.uni-koeln.de}}
  \and Michel GRABISCH\thanks{Corresponding author. Paris School of Economics,
    University of Paris I, 106-112, Bd. de l'H\^opital, 75013 Paris,
    France. Tel. (33) 144-07-82-85, Fax (33)-144-07-83-01. Email:
    \texttt{michel.grabisch@univ-paris1.fr}. The corresponding author thanks the
    Agence Nationale de la Recherche for financial support under contract
    ANR-13-BSHS1-0010.}  }

\date{\today}

\maketitle

\begin{abstract}
Many important values for cooperative games are known to arise from least square optimization problems. The present investigation develops an optimization framework to explain and clarify this phenomenon in a general setting. The main result shows that \emph{every} linear value results from some least square approximation problem and that, conversely, \emph{every} least square approximation problem with linear constraints yields a linear value.

\medskip
This approach includes and extends previous results on so-called least square values and semivalues in the literature. In particular, is it demonstrated how known explicit formulas for solutions under additional assumptions easily follow from the general results presented here.
\end{abstract}



\newpage

\section{Introduction}
Approximation of high-dimensional quantities or complicated
functions by simpler functions with linear properties from low-dimensional spaces has countless applications in physics, economics, operations research {\it etc.} In these applications, the quality of the approximation is usually measured by the Gaussian principle of least squared error,  which is also the guiding optimality criterion in the present investigation. Our study addresses a particular
case of such an approximation context with many applications in different
fields related to operations research, namely decision theory, game theory and
the theory of pseudo-Boolean functions.

\medskip
Where $N$ is a finite set with $n =|N|$ elements and collection $2^N$ of subsets, a set function $v:2^N\rightarrow \R$ assigns to
every subset of $N$ a real number, and is by definition of exponential
complexity (in $n$). Identifying subsets of $N$ with their characteristic (incidence) vectors (and thus $2^N$ with $\{0,1\}^n$), a set function can be viewed as a so-called pseudo-Boolean function  $f:\{0,1\}^n\rightarrow \R$ (\emph{cf.} Hammer and Rudeanu  \cite{haru68}). Of particular interest are those set functions which vanish on
the empty set, since they represent cooperative TU games with $N$ being
the set of players and the quantities $v(S)$ expressing the
benefit created by the cooperation of the members of $S\subseteq N$ (see, \emph{e.g.}, Peleg and Sudh\"olter \cite{pesu03}).  Under the additional stipulation of monotonicity, \emph{i.e.}, the property that $v(S)\leq v(T)$ holds whenever $S\subseteq T$, one
arrives at so-called capacities, which are a fundamental tool in the analysis of decision making under uncertainty
(\emph{cf.} Schmeidler~\cite{sch89}) or relative to several criteria (Grabisch and Labreuche~\cite{grla03e}).

\medskip
Being of exponential complexity, a natural question is to try to approximate
general set functions by simpler functions, the simplest being the additive set
functions, which are completely determined by the value they take on the $n$ singleton sets $\{i\}$ and are thus of linear complexity (in $n$).  In the field of pseudo-Boolean functions, the question has been addressed by Hammer and Holzman \cite{haho87} with respect to linear and
quadratic approximations, while approximation of degree $k$ was studied by
Grabisch \emph{et al.} \cite{grmaro99a}. In decision theory, linear
approximation amounts to the approximation of a capacity $\mu$ by a probability
measure $P$ (an additive capacity satisfying the additional constraint that $P(N)=1$).

\medskip
In game theory, the approximation of a game $v$ by an additive game (equivalently by a (payoff) vector in $\R^N$) is related to the concept of value or solution of a game: given $v$, find $x\in \R^N$ such that $\sum_{i\in N}x_i=v(N)$ and the $x_i$ represent as faithfully as
possible the contribution of the individual players $i$ in the total benefit $v(N)$. A very
natural approach for a value is to define it as the best least square
approximation of $v$, under the constraint $\sum_{i\in N}x_i=v(N)$, the
approximation being possibly weighted. Such values are called least square
values. An early and important contribution to this cooperative solution concept is due to Charnes {\it et al.} \cite{chgokero88}, who gave the general solution for the
weighted approximation with nonnegative weights, and exhibited the well-known Shapley
value \cite{sha53} as a least square value. Ruiz \emph{et al.} \cite{ruvaza98}, for example, generalized this approach and derived further values from least square approximation.

\medskip
The aim of this paper is to provide a general view on the set function approximation problem by placing it in the context of quadratic optimization and bringing well-known tools of convex analysis to bear on the problem. This approach not only generalizes existing results but also points to interesting connections and facts. Our formulation will remain general, although we will adopt most of the time the notation and ideas from cooperative game theory, due
to the great interest in this field towards values and how to obtain them.

\medskip
Our main result exhibits, roughly speaking, linear values and least square values for cooperative games to represent two sides of the same coin: we find that \emph{every} least square problem under linear constraints yields a linear value and that \emph{every} linear value arises as such a least square value (Section~\ref{sec:Coop-values}).

\medskip
The paper is organized as follows. Section~\ref{sec:ls-approximations} describes 
the general problem of least square approximation and gives the fundamental
result which will be used in the sequel (Theorem~\ref{t.main-theorem}). Section~\ref{sec:ls-values} concentrates on least square values, and establishes explicit solution formulas under mild conditions on the weights
used in the approximation. This model generalizes the approach to the Shapley value and to an optimization problem given in Ruiz \emph{et al.} \cite{ruvaza98}. We remark that, interestingly, the weights do not necessarily have to be all positive in our model.  Finally, we show in Section~\ref{sec:prob-values} how Weber's~\cite{Weber88} so-called
probabilistic values arise naturally in the present context.

\section{Least square approximations and linear operators}\label{sec:ls-approximations}
 We begin by reviewing some basic facts from convex optimization\footnote{see, \emph{e.g.}, Faigle \emph{et al.}  \cite{fakest02} or any other textbook for more details}. For integers $k,m\geq 1$, we denote by $\R^k$ the vector space of all $k$-dimensional (column) vectors and by $\R^{m\times k}$ the vector space of all ($m\times k$)-matrices $M=[m_{ij}]$ with coefficients $m_{ij}$. Generally, $M^T$ denotes the transpose of a matrix (or coefficient vector) $M$.

\medskip
Recall that any positive definite ($k\times k$)-matrix $Q =[q_{ij}]$ defines an inner product  {\it via}
$$
\(x|y\)_Q = x^TQy = \sum_{j=1}^k\sum_{i=1}^k q_{ij} x_iy_j 
$$
with the associated \emph{$Q$-norm} $\|x\|_Q = \sqrt{\(x|x\)_Q}$ on $\R^k$. Note that the choice $Q=I$ of the identity matrix $I$ yields the usual euclidian norm $\|x\|= \|x\|_Q$.

\medskip
Fix now a matrix $A^{m\times k}$, a linear map $b:\R^k\to \R^m$ as well as a linear map $c:\R^k\to \R^k$. For any $v\in \R^k$, denote by $\hat{v} =\hat{v}(A,b,c)$ the optimal solution of the quadratic minimization problem
\begin{equation}\label{eq.c-LSA}
\min_{Ax=b(v)} \|c(v)-x\|_Q^2.
\end{equation}
So, if the system $Ax=b(v)$ of linear equations has at least one solution, $\hat{v}$ is the (uniquely determined) best approximation of $c(v)$ in the solution space of $Ax=b(v)$ in the norm $\|\cdot\|_Q$. The key observation in our analysis is:

\medskip
\begin{lemma}\label{l.key-observation} Assume that $Ax=b(v)$ has a solution for every $v\in \R^k$ and that the map $c:\R^k\to\R^k$ is linear. Then $v\mapsto \hat{v}$ is a well-defined linear operator.
\end{lemma}

\Pf  Problem (\ref{eq.c-LSA}) is equivalent to the quadratic optimization problem
\begin{equation}\label{eq.c-LSA1}
\min_{Ax=b(v)} \frac12x^TQx - c(v)^Tx.
\end{equation}
Given that $Q$ is positive definite, it is well-known that $x$ is the unique optimal solution for problem (\ref{eq.c-LSA1}) if and only if there is a vector $y$ such that the associated Karush-Kuhn-Tucker (KKT) system
\begin{equation}\label{eq.KKT1}
\begin{array}{rcccl}
Qx & +&A^Ty &=& c(v)\\
Ax &&& = &b(v)
\end{array}
\end{equation}
is satisfied. Since $b$ and $c$ are linear functions in $v$, one immediately deduces from (\ref{eq.KKT1}) that also the optimal solutions of (\ref{eq.c-LSA}) are linear functions in $v$.

\qed

\medskip
\begin{theorem}\label{t.main-theorem} The operator $f:\R^k\to \R^k$ is linear if and only if there is a matrix $A\in \R^{m\times k}$, a linear function $b:\R^k\to\R^m$ and a linear function $c:\R^k\to \R^k$ such that $Ax=b(v)$ is always solvable and $f(v) = \hat{v}$ holds.
\end{theorem}

\Pf Lemma~\ref{l.key-observation} shows that the condition of the Theorem is sufficient for $f$ to be linear. Conversely, any $f(v)$ is obviously the optimal solution of the problem
$$
  \min_{x\in\R^k} \|f(v)-x\|^2_Q.
$$
Hence the Theorem is satisfied with the choice $A=0$ and $b=0$, for example.

\qed

\section{Values of cooperative games}\label{sec:Coop-values}
Let $N$ be a set of \emph{players} of finite cardinality $n=|N|$ and let $\cN$ be the collection of non-empty subsets $S\subseteq N$. A \emph{cooperative TU game} is a function $v:\cN\to \R$ (which is usually thought to be extended to all subsets of $N$ {\it via} $v(\emptyset) = 0$). So the set $\cG=\R^\cN$ of all cooperative TU games on $N$ is a vector space and isomorphic to $\R^k$ with $k=|\cN| = 2^n-1$.

\medskip
The \emph{additive (cooperative) games} correspond to those members $x\in\R^\cN$ that satisfy the homogeneous system of linear equations
$$
x(S) -\sum_{i\in S}x_i  = 0 \quad(S\in \cN)
$$
and one may be interested in the approximation of a game $v\in\cG$ by an additive game  with certain properties.  More general approximations might be of interest. For example, the linear constraints
\begin{eqnarray*}
\sum_{i\in N} x_i &=& v(N) \\
\sum_{S\in\cN}  x(S) &=&\sum_{S\in\cN} v(S)
\end{eqnarray*}
would stipulate an approximation of $v$ by a game that induces an efficient value (the first equality) and, furthermore, preserves the total sum of the $v(S)$ (second equality). Since the right-hand-side constraints are linear in $v$ Lemma~\ref{l.key-observation} says that least square approximations of this type are linear in $v$.

\medskip
A function $\Phi:\cG\to \R^N$ is \emph{value}  for $\cG$. It is straightforward, to view $\Phi(v)$ actually as an additive game that assigns the worth $\Phi_S(v)$ to the set $S$ of players by setting
$$
    \Phi_S(v) = \sum_{i\in S} \Phi_i(v) \quad(S\in \cN).
$$
Conversely, every additive game $v$ arises from  a parameter vector $\varphi\in \R^N$ so that
$$
     v(S) = \sum_{i\in S} \varphi_i \quad(S\in \cN).
$$
Hence the space of additive games is isomorphic with $\R^N$. Consequently, Theorem~\ref{t.main-theorem} implies that the \emph{linear} values are those which arise from  least square approximation problems with linear constraints.

\section{Least square values}\label{sec:ls-values}
We have seen that every linear value $\Phi:\cG\to\R^N$ can be interpreted as arising from a least square approximation problem. Special cases of seemingly more general least square problems have received considerable attention in the literature and led to the concept of \emph{least square values} and \emph{semivalues}. Take, for example, the weighted least square problem
\begin{equation}\label{eq:13}
\min_{x\in \R^N}\sum_{S\in \cN}\alpha_S(v(S) - x(S))^2 \quad\text{ s.t. }\quad \sum_{i\in N} x_i =v(N),
\end{equation}
where we set $x(S)=\sum_{i\in S} x_i$. So (\ref{eq:13}) asks for the best ($\alpha$-weighted) least square approximation of a game $v$ by an additive game $x$ under the additional efficiency constraint $x(N)=v(N)$.

\medskip
This problem has a long history. Hammer and Holzman (\cite{haho87})\footnote{later published in \cite{haho92}} studied both the above version and the unconstrained version with equal weights ($\alpha_S=1$ $\forall
S$), and proved that the optimal solutions of the unconstrained version yield
the Banzhaf value \cite{banzhaf65} (see also Section~\ref{sec:prob-values} below). More general versions of the un\-con\-strained problem were solved by Grabisch {\it et al.} \cite{grmaro99a} with the approximation being relative to the space of $k$-additive games (\emph{i.e.}, games whose
M\"obius transform vanishes for subsets of size greater than $k$)\footnote {see also Ding \cite{dilachchma10}, and Marichal and Mathonet \cite{mama11a}}.

\medskip
In 1988, Charnes {\it et al.} \cite{chgokero88} gave a solution for the case with the coefficients $\alpha_S$ being \emph{uniform} (\emph{i.e.}, $\alpha_S=\alpha_T$ whenever $|S|=|T|$) and strictly positive. As a particular case, the Shapley value was shown to result from the coefficient choice
\begin{equation}\label{eq.Charnes}
\alpha_S = \alpha_s =\binom{n-2}{s-1} = \frac{(n-2)!}{(s-1)!(n-1-s)!} \qquad (s=|S|).
\end{equation}

\medskip
\Rem. \emph{Ruiz {\it et al.} \cite{ruvaza98} state that problem (\ref{eq:13})  has a unique optimal solution for any choice of weights (see Theorem~3 there). In this generality, however, the statement is not correct as neither the existence nor the uniqueness can be guaranteed. So additional assumptions on the weights must be made.}

\medskip
We will first present a general framework for dealing with such situations and then illustrate it with the example of regular weight approximations and probabilistic values.

\subsection{Weighted approximation}\label{sec:weighted-approx}
For the sake of generality, consider a general linear subspace $\cF\subseteq \R^\cN$
of dimension $k=\dim\cF$, relative to which the
  approximation will be made.

\medskip
Let $W=[w_{ST}]\in \R^{\cN\times \cN}$ be a given matrix of weights $w_{ST}$. Let $c:\R^\cN\to\R^{\cN}$ be a linear function and consider, for any game $v$, the optimization problem
\begin{equation}\label{eq.valuation-opt-problem}
\min_{u\in \cF}~ (v-u)W(v-u)^T + c(v-u)^T \quad\mbox{with $c=c(v)$},
\end{equation}
which is equivalent with
\begin{equation}\label{eq.valuation-opt-problem1}
\min_{u\in \cF}~ uWu^T -\tilde{c}u^T,
\end{equation}
where $\tilde{c}\in \R^\cN$ has the components $\tilde{c}_S = c_S+2\sum_T
w_{ST}v_T$. A further simpli\-fication is possible by choosing a basis
$B=\{b_1,\ldots,b_k\}$ for $\cF$. With the identification
$$
 x=(x_1,\ldots,x_k)\in \R^k \quad\longleftrightarrow\quad u=\sum_{i=1}^k x_ib_i \in \cF,
$$
problem (\ref{eq.valuation-opt-problem1}) becomes
\begin{equation}\label{eq.valuation-opt-problem2}
\min_{x\in\R^k}~ \sum_{i=1}^k\sum_{j=1}^k q_{ij}x_ix_j - \sum_{i=1}^k\ov{c}_ix_i
\end{equation}
with the coefficients
$$
  q_{ij} = \sum_S\sum_T w_{ST}b_i(S)b_j(T) \quad\mbox{and}\quad \ov{c}_i =\sum_S \tilde{c}_Sb_i(S).
$$
Note that $\ov{c}:{\R^\cN}\to \R^k$ is a linear function.

Let $A\in \R^{m\times k}$ {be } a constraint matrix and $b:{
  \R^\cN}\to \R^m$ a linear function such that $Ax=b(v)$ has a solution for
every $v\in {\R^\cN}$. {If $Q=[q_{ij}]\in \R^{k\times k}$ is positive
  definite}, the problem
\begin{equation}\label{eq.valuation-opt-problem3}
\min_{x\in\R^k}~ \sum_{i=1}^k\sum_{j=1}^k q_{ij}x_ix_j -
\sum_{i=1}^k\ov{c}_ix_i\quad\mbox{s.t.}\quad Ax =b
\end{equation}
has a unique optimal solution $x^*$ which is linear in $v$ (Lemma~\ref{l.key-observation}). So we obtain the linear value $v\mapsto \hat{v}$ with components
$$
\hat{v}_j = u^*_{\{j\}} \quad(j\in N)\;,\;\;\mbox{where}\quad  u^* = \sum_{i=1}^k x_i^*b_i \in \cF.
$$

\medskip
In the model (\ref{eq:13}), for example, $\cF$ is the space $\cC$ of all additive games and has dimension $n$. The matrix $W$ is diagonal with the diagonal elements $w_{SS} = \alpha_S$. If $\alpha_S >0$ holds for all $S$, then $W$ is positive definite and the linearity of the implied value $v\mapsto \hat{v}$ follows directly from
Lemma~\ref{l.key-observation}.

\medskip
Otherwise, let us choose for $B$ the basis of unanimity games $\zeta_i$, $i\in N$, for $\cC$, where
$$
     \zeta_i(S) = \left\{\begin{array}{cl} 1&\mbox{if $i\in S$},\\
     0 &\mbox{if $i\not\in S$.}\end{array}\right.
$$
The associated matrix $Q=[q_{ij}]$ in model (\ref{eq:13}) has the coefficients
\begin{equation}\label{eq.semi-values}
     q_{ij} = \sum_{S\in\cN} \alpha_S \zeta_i(S)\zeta_j(S) =\sum_{S\ni\{i,j\}} \alpha_S.
\end{equation}
 For establishing a linear value, it suffices that $Q$ be positive definite, which is possible even when some of the $\alpha_S$ are negative (see Examples \ref{ex.symm} and \ref{ex.symm-pd} below).

\subsection{Regular weights}\label{sec:regular-weights}
While Lemma~\ref{l.key-observation} guarantees the existence of linear values resulting from approximation, explicit formulas can be given under additional assumptions on the weights. Restricting ourselves to objectives of type
$$
   \sum_{S\in \cN} \alpha_S(v_S-u_S)^2 +\sum_{S\in \cN} c_Su_S,
$$
we propose a simple framework that nevertheless includes all the cases treated in the literature so far. We say that the weights $\alpha_S$ are \emph{regular} if the resulting matrix $Q$ has just two types of coefficients $q_{ij}$, \emph{i.e.}, if there are real numbers $ p,q$ such that
$$
q_{ij} =\left\{\begin{array}{cl} q &\mbox{if $i=j$}\\ p &\mbox{if $i\neq j$.}\end{array}\right.
$$

\medskip
\begin{ex}\label{ex.symm} Assume that the weights $\alpha_S$ are uniform and set $\alpha(|S|) =\alpha_S$ . Then formula (\ref{eq.semi-values}) yields
$$
q_{ij} =\sum_{s=2}^n \binom{n-2}{s-2}\alpha(s) \quad\mbox{and}\quad q_{ii} = \sum_{s=1}^n \binom{n-1}{s-1}\alpha(s)
$$
holds for all $i\neq j$. So $Q=[q_{ij}]$ is regular.
\end{ex}

\medskip
\begin{lemma} Let $Q=[q_{ii}]\in \R^{k\times k}$ be regular with $q=q_{ii}$ and $p=q_{ij}$ for $i\neq j$. Then $Q$ is positive definite if and only if $q> p\geq 0$.
\end{lemma}

\Pf
For any $x\in \R^k$, we have after some algebra
\[
x^TQx = (q-p)\sum_{i=1}^k x_i^2 + p\overline{x}^2
\]
where $\overline{x}=\sum_{i=1}^n x_i$, which makes the claim of the Lemma obvious.

\qed

\medskip
Note that our model allows for possibly negative uniform coefficients, as shown in the following example.
\begin{ex}\label{ex.symm-pd}
Let $n=3$. We get $p=\alpha_2+\alpha_3$ and
$q=\alpha_1+2\alpha_2+\alpha_3$. Letting $\alpha>0$, the following vectors
$(\alpha_1,\alpha_2,\alpha_3)$ lead to a positive definite matrix $Q$:
\[
 (0,\alpha,0), \quad (\alpha,0,\alpha),\quad (0,\alpha,-\alpha), \text{ etc.}
\]
\end{ex}

\medskip
 For the remainder of this section, let $Q\in\R^{N\times N}$ be a regular matrix with parameters $q> p\geq 0$, $c\in \R^N$ a vector and $g\in \R$ a scalar. Setting $\1^T=(1,1,\ldots,1)$, the optimization problem

\begin{equation}\label{eq.reg-opt}
\min_{x\in \R^N}~ x^TQx- c^Tx \quad\mbox{s.t.}\quad \1^T x  =x(N) =g
\end{equation}
has a unique optimal solution $x^*\in\R^N$. Moreover, there is a unique scalar $z^*\in\R$ such that $(x^*,z^*)$ is the unique solution of the associated KKT-system
\begin{equation}\label{eq.KKT-reg}
\begin{matrix} Qx &-&z\1 &=& c/2\\
                \1^Tx &&&= &g.\end{matrix}
\end{equation}

Verifying this KKT-system, the proof of the following explicit solution formulas is straightforward.

\medskip
\begin{theorem}\label{t.regular} If $Q$ is regular, the solution $(x^*,z^*)$ of the KKT-system (\ref{eq.KKT-reg}) is:
\begin{eqnarray*}
z^* &=& (2(q +(n-1)p)g -C)/n\quad\mbox{(with $C = c\1^T =\sum_{i\in N} c_i$)}\\
x^*_i &=& (c_i +z^* -2pg)/(2q-2p)\quad(i\in N).
\end{eqnarray*}
If $Q$ is furthermore positive definite, then  $x^*$ is an optimal solution for (\ref{eq.reg-opt}).

\qed
\end{theorem}

\medskip
In the case of uniform weights $\alpha(s)$, the formulas in Theorem~\ref{t.regular} yield the formulas derived by Charnes {\it et al.} \cite{chgokero88} for problem (\ref{eq:13}). To demonstrate the scope of Theorem~\ref{t.regular}, let us look at the extremal problem\footnote{see also Sun {\it et al.} \cite{SuHaXu13} for similar problems}  studied by
Ruiz {\it et al.} \cite{ruvaza98a}

\medskip
\begin{equation}\label{eq:18}
\min_{x\in\R^N}\sum_{S\subseteq N} m_Sd(x,S)^2 \text{ s.t. } x(N)=v(N),
\end{equation}
where $m_S>0$ and
\[
d(x,S) = \frac{v(S) -x(S)}{|S|} - \frac{v(N\setminus S) - x(N\setminus S)}{n-|S|}.
\]
Letting $v^*(S) = v(N)-v(N\setminus S)$ and
\[
\overline{v}(S) = \frac{(n-|S|)v(S) +|S|v^*(S)}{n}
\]
(and thus $n\overline{v}(N)=v(N)$), we find that problem (\ref{eq:18}) becomes
\[
\min_{x\in\R^N}\sum_{S\subseteq N}\alpha_S(\overline{v}(S) - x(S))^2\text{
    s.t. } x(N)=n\overline{v}(N).
\]
with $\alpha_S= n^2m_S(|S|^2(n-|S|)^2)^{-1}$.
Because $v\mapsto \overline{v}$ and $v\mapsto g(v) = n\ov{v}(N)$ are  linear mappings, the optimal solutions of (\ref{eq:18}) yield an efficient linear value for any choice of parameters $m_S$ such that the associated matrix $Q$ is positive definite.

\medskip
If furthermore the weights $m_S$ (and hence the $\alpha_S$) are uniform, $Q$ is regular and the optimal solution can be explicitly computed from the formulas of Theorem~\ref{t.regular}.

\section{Probabilistic values}\label{sec:prob-values}
Weber~\cite{Weber88} introduced the idea of a \emph{probabilistic} value arising as the expected marginal contribution of players relative to a
probability distribution on the coalitions. For example, a \emph{semivalue} is a
  probabilistic value relative to probabilities that are equal on coalitions
  of equal cardinality.

\medskip
For our purposes, it suffices to think of the
\emph{marginal contribution} of an element $i\in N$ as a linear functional
$\partial_i:\cG\to {\R}$, where  $\partial_i^v(S)$ is interpreted as the marginal contribution of $i\in N$ to the coalition $S\subseteq N$ relative to the characteristic function  $v$.

\medskip
Probabilistic values can be studied quite naturally in the context of weighted  approximations. Indeed, let $p$ be an arbitrary probability distribution on $\cN$. Then the expected marginal contribution of $i\in N$ relative to the game $v$ is
$$
E(\partial^v_i) = \sum_{S\subseteq N} \partial^v_i(S)p_S.
$$

Let $\mu_i\in\R$ be an estimate value for the marginal contribution of $i\in N$. Then the expected observed deviation from $\mu_i$ is
$$
    \sigma(\mu_i) =\sqrt{\sum_{S\in\cN}p_S(\partial_i^v(S)-\mu_i)^2}.
$$
A well-known fact in statistics says that the deviation function $\mu_i\mapsto \sigma(\mu_i)$ has the unique minimizer $\mu =E(\partial^v_i)$, which can also be immediately deduced from the KKT conditions for the least square problem
$$
    \min_{\mu\in \R}~ \sum_{S\in\cN}p_S(\partial_i^v(S)-\mu)^2.
$$

\medskip
\paragraph{The values of Shapley and Banzhaf.} Shapley's~\cite{sha53} model assumes that player $i$ contributes to a coalition $S$ only if $i\in S$ holds and that, in this case, $i$'s  marginal contribution is evaluated as
$$
    \partial_i^v(S) = v(S)-v(S\setminus i).
$$
So only coalitions in $\cN_i=\{S\subseteq N\mid i\in S\}$ need to be considered. In order to speak about the ''average marginal contribution'', the model furthermore assumes:
\begin{enumerate}
\item The cardinalities $|X|$ of the coalitions $X\in \cN_i$ are distributed uniformly.
\item The coalitions $X\in\cN_i$ of the same cardinality $|X| =s$ are distributed uniformly.
\end{enumerate}

Under these probabilistic assumptions, the coalition $S\in \cN_i$ of cardinality $|S|=s$ occurs with probability
\begin{equation}\label{eq.Shapley-probability}
     p_S = \frac1n\cdot\frac{1}{\binom{n-1}{s-1}} = \frac{(s-1)!(n-s)!}{n!},
\end{equation}
which exhibits the Shapley value as a probabilistic (and hence approximation) value:
$$
    \sum_{S\in \cN_i} p_S[v(S)-v(S\setminus i)] = \sum_{S\in \cN} p_S[v(S)-v(S\setminus i)] =\Phi^{\rm Sh}_i(v).
$$

\medskip
\Rem. Among the probabilistic values, the Shapley value can also be characterized as the one with the largest entropy (Faigle and Voss~\cite{FaigleVoss11}).

\medskip
In contrast to the Shapley model, the assumption that all coalitions in $\cN_i$ are equally likely assigns to any coalition $S\in\cN_i$ the probability
\begin{equation}\label{eq.Banzhaf-probability}
     \ov{p}_S = \frac{1}{2^{n-1}}
\end{equation}
with the Banzhaf value \cite{banzhaf65} as the associated probabilistic value:
$$
    \sum_{S\in \cN_i} \ov{p}_S[v(S)-v(S\setminus i)] = \sum_{S\in \cN} \ov{p}_S[v(S)-v(S\setminus i)] = B^v_i.
$$


\end{document}